\documentclass[12pt]{article}
\usepackage{theorem,amsfonts}

\newtheorem{theorem}{Theorem}
\newtheorem{proposition}{Proposition}
\newtheorem{lemma}{Lemma}
\newtheorem{corollary}{Corollary}
\theorembodyfont{\upshape}

\newtheorem{remark}{Remark}
\newtheorem{proof}{Proof}

\newtheorem{keywords}{KEYWORDS}

\newtheorem{acknowledgement}{Acknowledgement}

\newcommand{\bt}{\begin{theorem}}
\newcommand{\et}{\end{theorem}}
\newcommand{\bl}{\begin{lemma}}
\newcommand{\el}{\end{lemma}}
\newcommand{\bp}{\begin{proposition}}
\newcommand{\ep}{\end{proposition}}
\newcommand{\br}{\begin{remark}}
\newcommand{\er}{\end{remark}}
\newcommand{\bc}{\begin{corollary}}
\newcommand{\ec}{\end{corollary}}
\newcommand{\be}{\begin{enumerate}}
\newcommand{\ee}{\end{enumerate}}
\newcommand{\bo}{\begin{proof}}
\newcommand{\eo}{\end{proof}}

\title{On heredity of strongly proximal actions}

\author{C. R. E. Raja}

\date{}

\begin{document}
\maketitle

\let\epsi=\epsilon
\let\vepsi=\varepsilon
\let\lam=\lambda
\let\Lam=\Lambda 
\let\ap=\alpha
\let\vp=\varphi
\let\ra=\rightarrow
\let\Ra=\Rightarrow 
\let\Llra=\Longleftrightarrow
\let\Lla=\Longleftarrow
\let\lra=\longrightarrow
\let\Lra=\Longrightarrow
\let\ba=\beta
\let\ga=\gamma
\let\Ga=\Gamma
\let\un=\upsilon

\begin{abstract}
We prove that action of a semigroup $T$ on compact metric space $X$ by 
continuous selfmaps is strongly proximal if and only if $T$ action on 
${\cal P}(X)$ is strongly proximal.  As a consequence we prove that affine
actions on certain compact convex subsets of finite-dimensional vector
spaces are strongly proximal if and only if the action is proximal.
\end{abstract}

{\footnotesize
\begin{keywords}
Proximal and strongly proximal actions, and probability measures.
\end{keywords}
}

Let $X$ be a complete separable metric space.  Let $T$ be a semigroup 
acting on $X$ by continuous selfmaps.  
A system $(X,T)$ is a pair consisting of a complete separable metric 
space $X$ 
and a semigroup $T$ acting on $X$ by continuous selfmaps.  In such a 
situation 
$X$ is a called a $T$-space.

Two points $x$ and $y$ in a $T$-space $X$ are said to be {\it 
proximal} if there exists a sequence $(t_n)$ in $T$ such that 
${\rm lim}t_nx = {\rm lim}t_ny$.  

We say that a system $(X, T)$ is 
{\it proximal} or the action of $T$ on $X$ is {\it proximal} if any 
two points 
$x$ and $y$ in $X$ are proximal.  

It is easy to see that group of special linear automorphisms on 
$\mathbb R^n $ 
action on $\mathbb R^n$ is proximal and the compact group actions are 
not proximal.  

Let ${\cal P}(X)$ be the space of all regular Borel probability 
measures on 
$X$, equipped with the weak* topology with respect to all continuous 
bounded 
functions.  It may be seen that ${\cal P}(X)$ equipped with the weak* 
topology 
is a complete separable metric space (see [P]).  The map $x \mapsto 
\delta _x$, 
maps $X$ homeomorphically onto a closed subset $\delta _X$, 
say of ${\cal P}(X)$ (see [P]) where 
$\delta _x$ is the measure concentrated at the point $x$.  Suppose a 
semigroup 
$T$ acts on $X$ by continuous selfmaps.  Then the action of $T$ on 
$X$ extends 
to an action on ${\cal P}(X)$ in the following natural way, 
for any $\lam \in {\cal P}(X)$ and any $t \in T$, 
$t\lam (E) = \lam (t^{-1} E)$ for any Borel subset $E$ of $X$. 

We say that a system $(X,T)$ is {\it strongly proximal} or the 
action of $T$ on $X$ is {\it strongly proximal} if for any 
$\lam \in {\cal P}(X)$, there exists a sequence $(t_n) \subset T$ 
such that $t_n \lam \ra \delta _x$ for some $x \in X$.

By considering ${1\over 2} (\delta _x +\delta _y)$ for any $x ,y \in X$, 
it is easy to see that any strongly proximal system is proximal; 
see [G] for more details on proximal and strongly 
proximal systems.  
But not all proximal systems are strongly proximal.  The action of 
the special Linear group $SL(V)$ on $V$ is proximal but it is not 
strongly proximal.

We now prove the following interesting result which is needed in the 
sequel.

\bp\label{sdd:prop:a1}
Let $T$ be a semigroup acting on a complete separable metric space 
$X$ by continuous selfmaps.  Then the action of $T$ on $X$ is 
strongly proximal 
if and only if the action of $T$ on ${\cal P}(X)$ is proximal.
\ep

\bo
Suppose the action of $T$ on $X$ is strongly proximal.  Let $\lam _1$ 
and $\lam _2$ be in ${\cal P}(X)$.  Then there exists a sequence 
$(t_n)$ in 
$T$ such that $t _n ({1\over 2}(\lam _1 + \lam _2) ) \ra \delta _x$ 
for 
some $x \in X$.  Then for given $1> \epsi >0$ there exists a compact 
subset 
$K$ of $X$ such that $$t_n \lam _1(K) + t_n\lam _2(K) > 2-\epsi \eqno 
(i)$$ 
for all $n \geq 1$.  Suppose for some $i =1,2$ and for some $m \geq 
1$, 
$t _{m} \lam _i (K) \leq 1-\epsi $.  Then since $\lam _1$ and $\lam 
_2$ are 
probability measures, we get that $$t_m\lam _1(K) + 
t_m\lam _2(K) \leq 2-\epsi$$ for some $m \geq 1$ which is a 
contradiction to 
$(i)$.  Thus, $$t_n \lam _i (K) > 1-\epsi$$ for $i=1,2$ and for all 
$n\geq 1$.  
By Prohorov's theorem (see [B] or [P]), the sequences $(t_n \lam _1)$ 
and 
$(t_n \lam _2)$ are 
relatively compact in ${\cal P}(X)$.  Let $\mu _1$ be a limit point 
of 
$(t_n \lam _1)$.  Then there exists a $\mu _2 \in {\cal P}(X)$ such 
that 
$${1\over 2}(\mu _1 + \mu _2) = \delta _x$$ and hence $\mu _1 =\delta 
_x$.  This 
implies that $${\rm lim} t_n \lam _1 = \delta _x = {\rm lim} t_n \lam 
_2 .$$  
Thus, the action of $T$ on ${\cal P}(X)$ is proximal.

Suppose the action of $T$ on ${\cal P}(X)$ is proximal.  Let 
$\lam \in {\cal P}(X)$.  Now for any $x \in X$, there exists a 
sequence $(t_n)$ 
in $T$ such that $${\rm lim} t_n \lam = {\rm lim} t_n\delta _x. \eqno 
(ii)$$  
For any $n \geq 1$, $t_n x\in \delta _X$ which is a closed 
$T$-invariant set 
and hence ${\rm lim} t_n x \in \delta _X$.  Thus, $(ii)$ implies that 
$$t_n \lam \ra \delta _ y$$ for some $y \in X$.
\hfill{$\fbox {}$}
\eo

Let $(X,T)$ be a dynamical system where $X$ is a complete separable 
metric 
space and $T$ be topological semigroup.  Let us now consider the 
map $\Psi \colon {\cal P}({\cal P}(X)) \ra {\cal P}(X)$ defined as  
$$\Psi (\rho ) = \int _{{\cal P}(X)} y d\rho (y) \in {\cal P}(X)$$ 
for any 
$\rho \in {\cal P}({\cal P}(X))$.  

We first establish the following properties of $\Psi$.

\bp\label{sdd:prop:c1}
Let $X$, $T$ and $\Psi$ be as above.  Then 
\be
\item $\Psi $ is a continuous $T$-equivariant map,

\item $\Psi (\delta _y ) = y $ for all $ y \in {\cal P}(X)$,

\item for $\rho \in {\cal P}({\cal P}(X))$, $\Psi (\rho ) = \delta 
_x$ for 
some $x \in X$ implies $\rho$ is a point mass concentrated at the 
point $x$. 

\item suppose $X$ is a semigroup, then $\Psi$ is a semigroup 
homomorphism.
\ee
\ep

\bo
Since $X$, ${\cal P}(X)$ and ${\cal P}({\cal P}(X))$ are all 
metrizable, it is 
enough to prove sequential continuity of $\Psi$.  Let $(\rho _n)$ be 
a 
sequence in ${\cal P}({\cal P}(X))$ such that $\rho _n \ra \rho \in 
{\cal P}({\cal P}(X))$.  Let $\nu _n = \Psi (\rho _n)$ for all $n$ 
and 
$\Psi (\rho ) = \nu$.  Let $f$ be a bounded continuous function on 
$X$.  Then 
the function $y \mapsto y(f)$ is a continuous bounded function on 
${\cal P}(X)$ 
and hence since $\rho _n \ra \rho$ in ${\cal P}({\cal P}(X))$, we 
have 
$$\nu _n (f) = \int y(f) d\rho _n (y) \ra \int y(f) d\rho (y) = \nu 
(f).$$ 
Thus, $\nu _n \ra \nu$ in ${\cal P}(X)$.  This proves that $\Psi $ is 
continuous.  Since the action of $T$ on $X$ is by continuous maps, we 
have 
$$t \nu = \int t y d\rho (y) = \int y d(t \rho)(y) ,$$ that is 
$t\Psi (\rho ) = \Psi (t \rho )$.  Thus, verifying property (1) of 
the map $\Psi$.  
It is easy to verify property (2) of the map $\Psi$.

Suppose for $\rho \in {\cal P}({\cal P}(X))$, $\Psi (\rho )=\delta 
_x=\nu$, say 
for some $x \in X$.  Then for any $\epsi >0 $ and any bounded 
continuous 
function $f$ on $X$ such that $f\geq 0$ and $f(x)=0$, let 
$$\sigma (f, \epsi) = \{y \in {\cal P}(X) \mid y(f) >\epsi\}.$$  Then 
$$0=\nu (f) \geq \int _{y\in \sigma (f, \epsi)} y(f) d\rho (y) \geq 
\epsi 
\rho (\sigma (f, \epsi ))$$ and hence 
$\rho (\sigma (f, \epsi )) =0$.  It is easy to see that $\sigma (f, 
\epsi)$ 
is an open set for all continuous bounded $f$ and $\epsi >0$.  Now 
let $W$ be 
the set of all nonnegative bounded continuous functions $f$ on $X$ 
which vanish at $x$ and 
$$B = \cup _{f\in W} \cup _{n=1}^\infty \sigma (f, {1\over n}).$$  
Then $B$ is an open set.  

We now claim that $B\cup {x} = {\cal P} (X)$ and $x\not \in B$.  Let 
$\lam (\not = \delta _x) \in {\cal P}(X)$.  Then choose a compact set 
$K$ 
such that $\lam (K) >{1\over n}$ for some $n$ and $x \not \in K$.  
Since $x \not \in K$, 
there exists a continuous function $f$ on $X$ such that $0 \leq f 
\leq 1$, 
$f(x) =0 $ and $f(y) =1$ for all $y \in K$.  Then $\lam (f) > 
{1\over n}$ 
and hence $\lam \in \sigma (f, {1\over n}) \subset B$.  Thus, 
${\cal P}(X) = B\cup {x}$ and it is easy to see that $x \not \in B$.  

We now claim that $\rho (B) =0$.   Let $K$ be any compact set 
contained $B$.  
Then there exists a finite number nonnegative continuous function 
$f_1, f_2, 
\cdots , f_k$ and a finite set of integers $n_1, n_2, \cdots , n_k$ 
such that 
$K \subset \cup \sigma (f_i, {1\over n_i})$.  Since 
$\rho (\sigma (f, \epsi) )= 0$ for all $f \in W$ and $\epsi >0$, we 
have 
$\rho (K) = 0$ and hence since $K$ is any arbitrary compact subset 
contained in 
$B$, we have $\rho (B) =0$.  Thus, $\rho$ is the mass concentrated at 
the 
point $x \in X$.  Thus, verifying property (3) of the map $\Psi$.

Suppose $X$ is a semigroup.  The $\Psi$ is a semigroup homomorphism 
follows 
from the facts 
\be
\item $\Psi$ is affine, that is for any $0\leq \ap_1, \ap _2 , \cdots 
, 
\ap _n \leq 1$ with $\sum \ap_i =1$ and for $\rho _1, \rho _2 , 
\cdots , 
\rho _n \in {\cal P}({\cal P}(X))$, $\Psi (\sum \ap _i \rho _i) = 
\sum \ap _i 
\Psi (\rho _i)$, 

\item the set of measures with finite supports in ${\cal P}(X)$ is 
dense in 
${\cal P}({\cal P}(X))$ and 

\item $\Psi$ is continuous.
\ee
\hfill{$\fbox{}$}
\eo

We now prove that strongly proximal actions on compact metric spaces 
is hereditary in the following sense.

\bt\label{sdd:thm:c1}
Let $X$ be a compact metric space and $T$ be a topological semigroup 
acting on $X$.  Then the following are equivalent: 
\be
\item $T$ action on $X$ is strongly proximal;

\item the action of $T$ on ${\cal P}(X)$ is strongly proximal.
\ee
\et

\bo 
Suppose $(X,T)$ is strongly proximal.  Let $\rho \in {\cal P}({\cal 
P}(X))$, 
and let $\nu = \Psi (\rho)\in {\cal P}(X)$.  Then since the action of 
$T$ on $X$ is strongly proximal, there exists a sequence 
$(t_n)$ in $T$ such that 
$t_n \nu \ra \delta _x$, for some $x \in X$.  Since for each $n$, 
$t_n \rho \in {\cal P} ({\cal P}(X))$ which is a compact metrizable 
space, the 
sequence $(t_n \rho)$ is a relatively compact sequence.  Now let 
$\rho _o \in 
{\cal P}({\cal P}(X))$ be a limit point of $(t_n \rho )$.  Since 
${\cal P}({\cal P}(X))$ is a metrizable space, there exists a 
subsequence 
$(t_{k_n})$ of $(t_n)$ such that $t_{k_n} \rho \ra \rho _0$ in 
${\cal P}({\cal P}(X))$.  Let $\nu _0 = \Psi (\rho _0) \in {\cal 
P}(X)$.  Then 
by Proposition \ref{sdd:prop:c1}, 
$$t_{k_n} \nu = t_{k_n} \Psi (\rho ) = \Psi (t_{k_n}\rho ) \ra \Psi 
(\rho _0) 
= \nu_0$$ in ${\cal P}(X)$ and hence 
since $t_n \nu \ra \delta _x$ in ${\cal P}(X)$ which is a metric 
space, we get 
that $\nu_0 = \delta _x$.  Again by Proposition \ref{sdd:prop:c1}, 
$\rho _0$ is 
the mass concentrated at $x \in X$.  Thus, the relatively compact 
sequence 
$(t_n \rho )$ has a unique limit point and hence it converges to the 
point 
mass concentrated at $x \in X$.  This proves that (1) implies (2).  
That (2) 
implies (1) follows from Proposition \ref{sdd:prop:a1} and from the 
remark that 
strongly proximal actions are proximal.
\hfill{$\fbox{}$}
\eo

As a consequence we have the following corollary for certain affine
actions: {\it affine action} of a semigroup $T$ on a closed convex 
subset $X$ of a locally 
convex vector space is an action of $T$ on $X$ such that 
$t(ax+(1-a)y) = at(x)+(1-a)t(y)$ 
for all $x, y\in X$ and all $0 \leq a \leq 1$.  

\bc
Let $\{v_1, v_2 , \cdots , v_n\}$ be a set of linearly independent vectors
on a locally convex vector space $V$ over reals.  Let $X$ be the
convex hull of 
$\{v_1, v_2, \cdots ,v_n \} $.  Let $T$ be a 
topological semigroup acting on $X$ by affine surjective 
maps.  Then the $T$ action on $X$ is proximal if and only if the $T$
action on $X$ is strongly proximal.
\ec

\bo
Let $F=\{v_i \mid 1\leq i\leq n\}$.  Since $F$ is compact, it is easy to
see that $X$ is compact.  Since $T$ is an affine action of $X$ by 
surjective maps and $F$ is
the set of extreme points of $X$, $F$ is a $T$-invariant set.  Thus, $T$
acts on $F$ also.  

Let $f\colon {\cal P}(F)\ra X$ be defined by $$f(\lam )=\sum 
\lam (v_i)v_i$$ 
for all $\lam \in {\cal P}(F)$.  Suppose $(\lam _n)$ is a
sequence in ${\cal P}(F)$ such that $\lam _n \ra \lam $ in ${\cal P}(F)$.
Then $\lam _n(v_i) \ra \lam (v_i)$ for all $i$.  This implies that
$f(\lam _n )= \sum \lam _n(v_i) v_i \ra \sum \lam (v_i) v_i$.  Thus, $f$
is continuous.  

Given any point $x$ in $X$ there exist $0\leq \lam _i \leq 1$ for $i =1,
2, \cdots n$ such that $ x=\sum \lam _i v_i$ and $\sum \lam _i =1$.  Since
$\{v_i \mid 1\leq i \leq n \}$ is a linearly independent set $\lam _i$'s
are unique.  This implies that $f$ is a bijection.  Since ${\cal P}(F)$ is
compact, $f$ is a homeomorphism.  Since the action is affine it is easy to
verify that $f(t\lam) = tf(\lam )$ for all
$t \in T$.  Suppose the $T$ action on $X$ is proximal.  Then the $T$
action on ${\cal P}(F)$ is proximal.  Now by Proposition
\ref{sdd:prop:a1}, $T$ action on $F$ is strongly proximal and hence 
by Theorem \ref{sdd:thm:c1}, $T$ action on ${\cal P}(F)$ is strongly
proximal.  Since $f$ is a homeomorphism preserving the $T$-action, $T$
action on $X$ is strongly proximal.
\hfill{$\fbox{}$}
\eo

\br
It should be noted that the conclusion of Theorem \ref{sdd:thm:c1} is 
valid for any Polish space if the map $\Psi$ is proper.
\er

In general for a complete separable metric space $X$, it is not clear 
that $({\cal P}(X), T)$ is strongly proximal if $(X,T)$ is strongly 
proximal.  
However the system $({\cal P}(X), T)$ does not admit a non-trivial 
$T$-invariant measure in the following sense:

\bc\label{sdd:cor:c1}
Let $T$ be a topological semigroup acting strongly proximally on a 
complete 
separable metric space $X$.  Suppose $\rho \in {\cal P}({\cal P}(X))$ 
is 
$T$-invariant.  Then $\rho$ is the mass concentrated at a point $x 
\in X$.
\ec

\begin{acknowledgement}
The work was done when the author was a  post-doctoral fellow of the 
Regional Commission CCRRDT du Pays de Loire at Universite d'Angers,
France.  I would like to thank Dr. Piotr Graczyk for his kind hospitality
during my stay at Angers.
\end{acknowledgement}

\noindent {C. Robinson Edward Raja, \\ 
Stat-Math Unit,\\ 
Indian Statistical Institute,\\ 
8th Mile Mysore Road,\\
R. V. College Post,\\
Bangalore - 560 059.\\
India.  \\
e-mail: creraja@isibang.ac.in }


\begin{thebibliography}{F2}

\footnotesize

\bibitem [B]{B} P. Billingley, Convergence of Probability Measures, 
John Willey and Sons, New York-Tornoto, 1968.

\bibitem [G]{G} S. Glasner, Proximal flows on Lie groups, Israel 
Journal of Mathematics, 45 (1983), 97-99.

\bibitem [P]{P} K. R. Parthasarathy, Probability Measures on Metric 
Spaces, 
Academic Press, New York-London, 1967.
\end{thebibliography}
\end{document}